\documentclass[12pt]{amsart}
\usepackage{amssymb}

\begin{document}
\newcommand{\OP}{{\mathcal O}_{{\mathbb P}^n}}
\newcommand{\HH}{{\mathrm H}}
\newcommand{\OX}{{\mathcal O}_X}
\newcommand{\OZ}{{\mathcal O}_Z}
\newcommand{\ol}[1]{\overline{#1}}
\newcommand{\sta}[1]{\stackrel{#1}{\to}}
\newcommand{\ul}[1]{\underline{#1}}
\newtheorem{thm}{Theorem}[section]              
\newtheorem{propose}[thm]{Proposition}
\newtheorem{lemma}[thm]{Lemma}
\newtheorem{cor}[thm]{Corollary}
\newenvironment{rmks}{
			\trivlist \item[\hskip \labelsep{\it Remarks}:]
		      }{
		     	\endtrivlist}

\title[Arithmetically Cohen-Macaulay bundles]{Arithmetically Cohen-Macaulay bundles on hypersurfaces}
\author{N.~Mohan Kumar}
\address{Department of Mathematics, Washington University in St. Louis,
St. Louis, Missouri, 63130}
\email{kumar@wustl.edu}
\urladdr{http://www.math.wustl.edu/\~{} kumar}
\author{A.~P.~Rao}
\address{Department of Mathematics, University of Missouri-St. Louis,
St. Louis, Missouri 63121}
\email{rao@arch.umsl.edu}
\author{G.~V.~Ravindra}
\address{Department of Mathematics, Washington University in St. Louis,
St. Louis, Missouri, 63130}
\email{ravindra@math.wustl.edu}
\subjclass{14F05}
\keywords{Vector bundles, hypersurfaces, Arithmetically Cohen-Macaulay}
\begin{abstract}
We prove that any rank two arithmetically Cohen-Macaulay vector bundle
on a general hypersurface of degree at least three in $\mathbb{P}^5$
must be split.
\end{abstract}
\maketitle

\begin{section}{Introduction}
Let $X\subset\mathbb{P}^n$ be a hypersurface of degree $d$. A vector
bundle $E$ on $X$ is called {\it Arithmetically Cohen-Macaulay} (ACM
for short) if $\HH^i(E(k))=0$ for all $k$ and $0<i<n-1$. By Horrock's
criterion \cite{Horrocks}, this is equivalent to saying that $E$ has a
resolution,
$$0\to F_1\to F_0\to E\to 0,$$ where the $F_i$'s are direct sums of
line bundles on $\mathbb{P}^n$. If $d=1$, $E$ is a direct sum of line
bundles ({\it op. cit}); the ACM condition is vacuous for $n=1,2$.

In this article, we will be interested in ACM bundles of rank two. For
$n=3$, ACM rank two bundles are ubiquitous [Remark \ref{n3}]. Hence we
will deal with smooth hypersurfaces $X$ of degree $d\geq 2$ in
$\mathbb{P}^n$ with $n\geq 4$ and ACM rank two bundles on $X$ which
are not direct sums of line bundles of the form $\OX(k)$. By the
Grothendieck-Lefschetz theorem, these bundles are the same as the
indecomposable rank two ACM bundles on $X$.

Our main theorem is 
\begin{thm}\label{main}
Let $X\subset \mathbb{P}^n$ be a smooth hypersurface of degree $d\geq
2$.
\begin{enumerate}
\item If $n\geq 6$, any ACM rank two bundle on $X$ is a direct sum of
line bundles (Kleppe \cite{Kleppe}).
\item 
\begin{enumerate}
\item If $n=5$ and $X$ is smooth, then for any ACM rank two bundle $E$
on $X$, $\HH^1(X, E^{\vee}\otimes E(k))=0~~\forall~~ k$. In
particular, $E$ is rigid.
\item If $n=5$ and $X$ is general (for a dense Zariski open subset of
  the parameter space of hypersurfaces) of degree $d\geq 3$, then any
  ACM bundle of rank two is a direct sum of line bundles.
\end{enumerate}
\item If $n=4$ and $X$ is general of degree $d\geq 5$, then for any
ACM rank two bundle $E$ on $X$, $\HH^i(X, E^{\vee}\otimes E)=0$ for
$i=1,2$. In particular, $E$ is rigid.
\end{enumerate}

\end{thm}

\begin{rmks}
\hfill

\begin{enumerate}
\item Part 1) of the theorem is known by the work of Kleppe
(\cite{Kleppe}, see Proposition 3.2). Our proof is different and
essentially falls out of some of the computations necessary for part
2) of the theorem.
\item\label{construction} If $n=5$ and $d\geq 2$, there certainly
exist special smooth hypersurfaces with non-split ACM bundles.  Here
is a fairly simple way to construct them: Let $f,g,h$ be a regular
sequence of homogeneous polynomials. Let $a,b,c$ be non-constant
homogeneous polynomials such that these six polynomials have no common
zero in $\mathbb{P}^n$ and such that $F=af+bg+ch$ is homogeneous of
degree $d\geq 2$. Let $X$ be defined by $F=0$. Then we have an exact
sequence,

\noindent $0\to E\to \OX(-\deg f)\oplus\OX(-\deg g)\oplus\OX(-\deg h)\to I\to
0,$

where $I$ is the ideal generated by $f,g,h$ in $\OX$, the map to $I$
is the obvious one and $E$ is the kernel. One easily checks that $E$
is an indecomposable ACM bundle on $X$ of rank two. Any smooth quadric
hypersurface in $\mathbb{P}^5$ has a Pl\"{u}cker equation $F=0$ and so
the above construction applies.

\item In the case of $n=4$, the rigidity statement was proved for
quintic threefolds by Chiantini and Madonna
\cite{madonna2001}. Further, it is known that indecomposable ACM
bundles of rank $2$ exist for any smooth hypersurface of degree $d$
with $2\leq d\leq 5$ (\cite{AR}, \cite{Beauville}, \cite{IM}). One way
to see this is to note that any such hypersurface contains a line and
hence the construction in remark (\ref{construction}) applies. On the
other hand, it was shown by Chiantini and Madonna \cite{madonna2004}
that such bundles do not exist for a general sextic in $\mathbb{P}^4$
and one expects the same to be true for general hypersurfaces of
degree $d\geq 6$ in which case our result is trivially true.

\item \label{n3} When $n=3$, any smooth hypersurface contains a point
and hence the construction in remark (\ref{construction}) applies in
this case too.

\end{enumerate}
\end{rmks}

We now give a brief sketch of the proof of (\ref{main}) (2b). Suppose
 we have a rank two indecomposable ACM bundle $E$ on a hypersurface
 $X$ of degree $d$.  We show that the module $N=\oplus \HH^{2}(X,
 E^{\vee}\otimes E(k))$ is a non-zero graded cyclic module generated
 by an element of degree $-d$ (Lemma \ref{monogenic}).  In the
 $\mathbb{P}^5$ case, $N$ is a Gorenstein module and the socle element
 is in degree $2d-6$. Thus we see that $N_k\neq 0$ for all $-d\leq
 k\leq 2d-6$.  By means of a deformation-theoretic argument, we will
 show that if $X$ is general, then the multiplication map $N_{-d}\to
 N_0$ by any $g\in\HH^0(X,\OX(d))$ is zero [Corollary
 \ref{vanishing1}], which by the cyclicity implies that $N_0=0$. This
 is a contradiction if $d$ is at least three.
\end{section}

\begin{section}{Cohomology computations}
We will work over an algebraically closed field of characteristic
zero, though most of the arguments will go through in characteristic
not equal to two.  We will assume throughout that $X$ is a smooth
hypersurface of degree $d\geq 2$ in $\mathbb{P}^n$ with $n\geq 4$. Let
$E$ be rank two bundle on $X$. By the Grothendieck-Lefschetz theorem,
$\mathrm{Pic}\,\mathbb{P}^n\to\mathrm{Pic}\, X$ is an isomorphism. So
by normalising $E$, we will assume that $c_1(E)=e$ where $e=0$ or
$-1$. We will now assume $E$ is ACM. Then we have a minimal
resolution
$$0\to F_1\to F_0\to E\to 0,$$ where the $F_i$'s are direct sums of
line bundles on $\mathbb{P}^n$ of rank $r$ ($E\cong
\OX(a)\oplus\OX(b)$ for some integers $a$ and $b$ if and only if
$r=2$).

Dualizing the above we
get,
$$0\to F_0^{\vee}{}\to F_1^{\vee}{}\to {\mathcal E}xt^1_{\OP}(E,\OP)\to 0.$$
Applying ${\mathcal H}om_{\OP}(E,*)$ to the exact sequence, 
$$0\to \OP\to\OP(d)\to\OX(d)\to 0,$$ we see that,
$${\mathcal E}xt^1_{\OP}(E,\OP)=E\,^{\vee}{}\,(d)=E(d-e).$$
This implies that $F_1^{\vee}{}=F_0(d-e)$. Thus we can rewrite our
minimal resolution as

\begin{equation}\label{basic1}
0\to F_1=F_0\,^{\vee}{}\,(e-d)\stackrel{\phi}{\to} F_0\to E\to 0
\end{equation}
In fact by \cite{Beauville}, the map $\phi$ can be chosen to be
skew-symmetric though we will not use this fact.

By restricting to $X$, we get an exact sequence,
\begin{equation}\label{basic2}
0\to E(-d)\to F_1\otimes\OX=\ol{F_1}\stackrel{\ol{\phi}}{\to}
F_0\otimes\OX=\ol{F_0}\to E\to 0 
\end{equation}
Let the image of $\ol{\phi}$ be denoted by $G$. Then $G$ is a vector bundle
of rank $r-2$ on $X$.

Taking exterior powers of $\phi$ in (\ref{basic1}), we have an exact
sequence
\begin{equation}\label{basic3}
0\to \Lambda^2 F_1\to \Lambda^2 F_0\to {\mathcal F}\to 0
\end{equation}
for some cokernel sheaf $\mathcal{F}$. 
\begin{lemma}\label{exterior}
We have an exact sequence
$$0\to L\to {\mathcal F}\to {\mathcal F}\otimes \OX=\ol{\mathcal F}\to
0$$ where $L$ is a line bundle on $X$ and $\ol{\mathcal F}$ is a
vector bundle of rank $2r-3$ on $X$ which fits into a natural exact
sequence
$$0\to E\otimes G\to \ol{\mathcal F}\to \OX(e)=\Lambda^2 E\to 0.$$

\end{lemma}
\begin{proof}
We certainly have an exact sequence
$$0\to {\mathcal I}{\mathcal F}\to {\mathcal F}\to \ol{\mathcal F}\to
0$$ where ${\mathcal I}$ is the ideal sheaf defining $X$ in
$\mathbb{P}^n$. It suffices to check that ${\mathcal I}{\mathcal F}$
is a line bundle on $X$ and $\ol{\mathcal F}$ is a vector bundle of
rank $2r-3$ for the first part of the lemma. Both statements are
local. Clearly ${\mathcal F}$ is set-theoretically supported only
along $X$. By localising we may assume that (\ref{basic1}) looks like,
$$0\to F_1\to F_0\to E\to 0$$ where the $F_i$'s are free of rank $r$
and the matrix of the map $F_1\to F_0$ is the diagonal matrix
$(x,x,1,1,\ldots,1)$ where $x=0$ defines $X$. Then the matrix in
(\ref{basic3}) is given by the diagonal matrix $(x^2,x,\ldots x,
1,1,\ldots,1)$ where we have one $x^2$, $2(r-2)$ $x$'s and the rest
1's. The claim about $\mathcal{I}\mathcal{F}$ and
$\ol{\mathcal{F}}$ follows easily from this.

To see the final exact sequence, we restrict (\ref{basic3}) to $X$ to
get

$$\Lambda^2\ol{F_1}\to\Lambda^2\ol{F_0}\to\ol{\mathcal F}\to 0.$$ 

From the exact sequence,
$$0\to G\to\ol{F_0}\to E\to 0$$ we note that

$$\mathrm{im}\,(\Lambda^2\ol{F_1}\to\Lambda^2\ol{F_0})=\mathrm{im}\,
(\Lambda^2G\to \Lambda^2\ol{F_0}).$$ 

This leads to the natural diagram,

\[
\begin{array}{ccccccccc}
&&0&&0&&&&\\
&&\downarrow&&\downarrow&&&&\\
&&\Lambda^2G&=&\Lambda^2G&&&&\\
&&\downarrow&&\downarrow&&&&\\
0&\to &{\mathcal G}&\to &\Lambda^2\ol{F_0}&\to &\Lambda^2E=\OX(e)&\to& 0\\
&&\downarrow&&\downarrow&&\|&&\\
0&\to&E\otimes G&\to&\ol{\mathcal F}&\to&\OX(e)&\to&0\\
&&\downarrow&&\downarrow&&&&\\
&&0&&0&&&& 
\end{array}
\]
where $\mathcal{G}$ is defined by the diagram.
\end{proof}
Next we note the vanishing of various cohomologies. From
(\ref{basic1}) and (\ref{basic3}) we have

\begin{equation}\label{basic4}
\HH^i(E(*))=\HH^i({\mathcal F}(*))=\HH^i(L(*))=0, 1\leq i\leq n-2
\end{equation}

This implies from the lemma above that

\begin{equation}\label{basic5}
\HH^i(\ol{\mathcal F}(*))=0, 1\leq i\leq n-3
\end{equation}

Tensoring the exact sequence, $0\to E(-d)\to\ol{F_1}\to G\to 0$ with
$E\,^{\vee}{}$ and taking cohomologies, we get, 

\begin{equation}\label{basic6}
\HH^i(E\,^{\vee}{}\otimes G(k))=\HH^{i+1}(E\,^{\vee}{}\otimes E(-d+k))
~~\forall~~ k,~~ 1\leq i\leq n-3
\end{equation}

Similarly, tensoring the exact sequence $0\to G\to \ol{F_0}\to E\to 0$
with $E\,^{\vee}{}$ and taking cohomologies, we get,
\begin{equation}\label{basic7}
\HH^i(E\,^{\vee}{}\otimes E(k))=\HH^{i+1}(E\,^{\vee}{}\otimes G(k))
~~\forall~~ k,~~ 1\leq i\leq n-3
\end{equation}

From the exact sequence in the lemma, using
(\ref{basic4}),(\ref{basic5}) we get,

\begin{equation}\label{basic8}
\HH^i(E\otimes G(k))=0, 2\leq i\leq n-3
\end{equation}

\begin{lemma}\label{suff1}
The vector bundle $E$ is a direct sum of line bundles if and only if
$\HH^2(E\, ^{\vee}{}\otimes E(-d))= 0.$
\end{lemma}

\begin{proof}
If $E$ is a direct sum of line bundles, this vanishing is clear. So
assume the vanishing.  From (\ref{basic6}) we see that
$\mathrm{Ext}\,^1(E,G)= \HH^1(E\,^{\vee}{} \otimes G)=0$ and thus we
see that the exact sequence $0\to G\to\ol{F_0}\to E\to 0$
splits. Since $\ol{F_0}$ is a direct sum of line bundles, we see that
so is $E$.
\end{proof}

\begin{cor}\label{monogenic}
If $E$ is an indecomposable bundle, then the finite length module
$\oplus_k 
\HH^1(E\otimes G(k))$ is a non-zero cyclic module generated by an
element of degree $-e$.
\end{cor}
\begin{proof}
By lemma \ref{exterior}, we have an exact sequence,
$$0\to E\otimes G\to \ol{\mathcal F}\to \OX(e)\to 0.$$ Taking
cohomologies and using (\ref{basic5}), we get the fact that the module
is cyclic generated by an element of degree $-e$. If it is zero then
$$\HH^1(E\otimes G(-e))=\HH^1(E^{\vee}\otimes G)=0$$ and so by 
(\ref{basic6}) $\HH^2(E^{\vee}\otimes E(-d))=0$. By the previous
lemma, the bundle would then have to be split.
\end{proof}
\end{section}

\begin{section}{Deformation criteria for acm vector bundles}
We have already noted that for $n\leq 5$, given any degree $d\geq 2$,
there exists a smooth hypersurface of degree $d$ and a non-split ACM
bundle of rank 2 on the hypersurface. So, in this section, we analyze
the situation, when there is such a vector bundle on a general
hypersurface $X$ of degree $d$.

We will start with some results on vector bundles on families of
varieties. We will not prove the most general results in this
direction, but just what we need. Most of the arguments are similar to
those used in the construction of quot schemes and are well-known to
experts. We are really interested in Corollary \ref{infinitesimal} and
much of what follows consists of technical results to achieve it.

Let us start by fixing some notation. All schemes will be of finite
type over the base field. Let $p:X\to S$ be a flat projective
morphism. If $q:T\to S$ is any morphism, we will denote by
$X_T=X\times_S T$, the fiber product and $p':X_T\to T$, $q':X_T\to X$,
the natural $S$-morphisms. We start with an elementary lemma, whose
proof we omit.

\begin{lemma}\label{openclosed}
Let $X,S,T$ be as above. Let $V\subset X$ be any subset such that
$q'(X_T)\cap V=\emptyset$. Then, $q'(X_T)\cap p^{-1}p(V)=\emptyset$.
\end{lemma}

\begin{propose}\label{spreading1}
Let $p:X\to S$ be a flat projective morphism. Let $F_1,F_2$ be two
vector bundles on $X$ with $\HH^1(X_s, F_2^{\vee}\otimes
F_1\mid_{X_s})=0~~ \forall ~~s\in S$. Let $m\geq 0$ be an integer.Then
there exists a scheme $q:S'\to S$ and an exact sequence,
$$q'^*F_2\to q'^*F_1\to G\to 0$$ where $G$ is a rank $m$ vector bundle
on $X_{S'}$. Furthermore, one has the following universal property:
for any reduced $S$-scheme $T$ with structure morphism $f:T\to S$ and
an exact sequence
$$f'^*F_2\to f'^*F_1\to G'\to 0$$ with $G'$ a rank $m$ vector bundle
on $X_T$, there exists an $S$-morphism $T\to S'$ such that the second
sequence is just the pull back of the first by the induced morphism
$X_T\to X_{S'}$.
\end{propose}
\begin{proof}
The hypothesis on $\HH^1$ ensures that
$\mathcal{E}=p_*(F_2^{\vee}\otimes F_1)$ is a vector bundle on
$S$. Let
$$\mathcal{H}=\mathbb{A}(\mathcal{E}^{\vee})\sta{q} S.$$
Then $q^*\mathcal{E}$ has a section and $\mathcal{H}$ is universal
with this property.
We have the fiber product diagram,
\[
\begin{array}{ccc}
X_{\mathcal{H}}&\sta{q'}&X\\
\downarrow p'&&\downarrow p\\
\mathcal{H}&\sta{q}&S
  \end{array}
\]

 By the flatness of $q$, we have,
$$p'_*q'^*(F_2^{\vee}\otimes F_1)=q^*\mathcal{E}$$ and so
$q'^*(F_2^{\vee}\otimes F_1)$ has a section. Thus on $X_{\mathcal{H}}$
we get the universal map $q'^*F_2\to q'^*F_1$. Let $G$ be the
cokernel. By semi-continuity, the points $x\in X_{\mathcal{H}}$ such
that $\dim_{k(x)}G\otimes k(x)<m$ constitute an open set, which we
denote by $V$. Since $p'$ is flat $p'(V)$ is open. Let
$\mathcal{H}'\subset\mathcal{H}$ be the closed subset with the reduced
scheme structure which is the complement of $p'(V)$. Thus, replacing
$\mathcal{H}$ with $\mathcal{H}'$, we get an exact sequence, where $G$
has the property that at every point $x\in X_{\mathcal{H}'}$,
$\dim_{k(x)} G\otimes k(x)\geq m$. Again, the set of points $x\in
X_{\mathcal{H}'}$ such that $\dim_{k(x)} G\otimes k(x)>m$ is a closed
subset, say $Z$.  Since $p'$ is proper, we may take
$S'\subset\mathcal{H}'$ to be the open set which is the complement of
the closed subset $p'(Z)$. It is clear that on $X_{S'}$ we have an
exact sequence as claimed.

We need to check the universal property. So, let $f:T\to S$ be as in
the proposition with the exact sequence as mentioned.  Let us denote
by $f':X_T\to X$ and $p'':X_T\to T$ the corresponding maps.  The
existence of the map $f'^*F_2\to f'^*F_1$ gives a section of
$p''_*f'^*(F_2^{\vee}\otimes F_1)$ which is equal to $f^*\mathcal{E}$
by semi-continuity. This implies that we have a morphism
$g:T\to\mathcal{H}$ over $S$ and the exact sequence
$$q''^*E_2\to q'^*E_1\to G\to 0$$
on $X_{\mathcal{H}}$ pulls back to the one on $X_T$ via the induced
map $g':X_T\to X_{\mathcal{H}}$. Since $G'$ is a rank $m$ vector
bundle, we see that $g'(X_T)\cap V=\emptyset$ and thus we see that by
lemma \ref{openclosed}, $g'(X_T)\cap p'^{-1}p'(V)=\emptyset$. Since
$T$ is reduced, this implies that $g$ factors through
$\mathcal{H}'$. Next, we notice that $g'(X_T)\cap Z=\emptyset$ and by
another application of lemma \ref{openclosed}, we are done.
\end{proof}

Given a rank two (non-split) ACM bundle $E$ on $X$, we rewrite
(\ref{basic1}) with $F_0\cong \oplus_{i=1}^{r}\OP (a_i)$, $a_1\geq
a_2\geq \cdots \geq a_{r}$ and $r>2$ to get

$$0 \to \oplus_{i=1}^r\OP(-a_i+e-d) \stackrel{\phi}\to
\oplus_{i=1}^r\OP(a_i) \to E \to 0$$

As before, we assume that $e=0$ or $e=-1$.

\begin{lemma}\label{finite}
When $n\geq 4$, for fixed $d$ there are only finitely many
possibilities for $\underline{a}:=(a_1,\cdots,a_r)$.
\end{lemma}

\begin{proof}
Since sequence (\ref{basic1}) is minimal, we see that all entries of
$\phi$ must be at least of degree 1 and thus $\deg\det\phi\geq r$. But
this determinant is just the square of the equation defining our
hypersurface and thus it must be $2d$. So $r\leq 2d$ and is
bounded. 
From the fact that $\phi$ is an inclusion, we see that $-a_r+e-d\leq
a_1$. Thus, if we show that $a_1$ is bounded above, then it will
follow that $a_r$ is bounded below and then we will have only finitely
many possibilities for the $a_i$'s.

For this we proceed as follows. Let $a=a_1$. Since (\ref{basic1}) is
minimal, we get an inclusion $\OX(a_1)\to E$ such that the quotient is
torsion free. Thus, we get an exact sequence,
$$0\to \OX(a)\to E\to I(e-a)\to 0$$ where $I$ is the ideal sheaf of a
codimension two subscheme $Z\subset X$. This implies
$\omega_Z\cong\OZ(e-2a+d-n-1)$ where $\omega_Z$ is the the canonical
bundle of $Z$. Let $\pi:Z\to \mathbb{P}^{n-3}$ be a general
projection, so that $\pi$ is finite. Then
$$\pi_*\omega_Z={\mathcal H}om(\pi_*\OZ,{\mathcal
O}_{\mathbb{P}^{n-3}}(-n+2)).$$ 

Since $\mathcal{O}_{\mathbb{P}^{n-3}}$ is a direct summand of
$\pi_{\ast}\OZ$, $\HH^0(\pi_*\omega_Z(n-2))\neq 0$ and thus
$\HH^0(\OZ(e-2a+d-3))\neq 0$. Since $n\geq 4$, and $E$ is ACM we see
that $\HH^1(I(*))=0$ and thus $\HH^0(\OX(k))\to \HH^0(\OZ(k))$ is
onto. Since $X$ is a smooth hypersurface, $\HH^0(\OX(k))=0$ for
$k<0$. Thus we see that $e-2a+d-3\geq 0$ or $a\leq (e+d-3)/2$. 
\end{proof}

\begin{thm}\label{spreading}
Let ${\mathcal P}$ denote the parameter space of all degree $d$ smooth
hypersurfaces in $\mathbb{P}^n$ with $n\geq 4$ and ${\mathcal
X}\subset \mathbb{P}^n\times \mathcal{P}$ the universal
hypersurface. Assume that there exist rank two ACM non-split bundles
on hypersurfaces corresponding to a Zariski dense subset of ${\mathcal
P}$. Then there exists a scheme ${\mathcal P}'$, a dominant morphism
${\mathcal P}'\to {\mathcal P}$ and a rank two bundle over
$\mathcal{X}\times_{\mathcal P}\mathcal{P}'$ which is ACM and
non-split for any point in $\mathcal{P}'$.
\end{thm}

\begin{proof}
By lemma \ref{finite}, for any rank two non-split ACM bundle on a
hypersurface $X$ of degree $d$, we have finitely many choices for
$r,e$ and the $a_i$'s. Let us fix one of these. Then for all vector
bundles $E$ with these invariants, we have a presentation,
$$F_0^{\vee}(e-d)\otimes \OX\to F_0\otimes \OX\to E\to 0.$$ Thus, we
consider the flat projective morphism $p:\mathcal{X}\to \mathcal{P}$
and the vector bundles $F_0^{\vee}(e-d), F_0$ pulled back to
$\mathcal{X}$, which we denote by the same name. Notice that the
vanishing condition on first cohomology in Proposition
\ref{spreading1} holds if $n\geq 3$.  Thus we get a scheme
$q:S'=S'(e,r,a_i)\to \mathcal{P}$ as in the proposition for $m=2$. Let
$G$ be the corresponding rank two bundle on $\mathcal{X}_{S'}$. We
have a closed subset $S''\subset S'$ where the map $F_0^{\vee}(e-d)\to
F_0$ is minimal. We may restrict to $S''$ and let
$p:\mathcal{X}_{S''}\to S''$ be the structure map.

By the relative version of Serre vanishing \cite{Mumford},there exists
an integer $m_0$ such that for all $m\geq m_0$ and all $i>0$,
$\mathrm{R}^ip_*G(m)=0$. Since $G$ is a vector bundle and $p$ is flat,
by repeated application of semi-continuity theorems [see for example,
page 41, \cite{Mumford}], one sees that $\HH^i(X_s,
G(m)_{\mid{X_s}})=0$ for all $i>0$ and all $m\geq m_0$. By duality
this is also true for all $i<n-1$ and all $m\leq m_1$ for a suitable
$m_1$.  Thus we see that there are only finitely many integers $k$
such that $\HH^i(X_s, G(k)\mid_{X_s})\neq 0$ for some $s\in S''$ and
some $i$ with $0<i<n-1$. Since the set of $s\in S''$ such that
$\HH^i(X_s,G(k)\mid_{X_s})\neq 0$ for fixed $i,k$ is a closed subset,
we see that there is a well-defined closed subset $Z\subset S''$ such
that $\HH^i(X_s,G(k)\mid_{X_s})\neq 0$ for some $k$ and some $i$ with
$0<i<n-1$ if and only if $s\in Z$. If we let $T=S''-Z$, we see that on
$\mathcal{X}_T$, the bundle $G$ has the property that on each fibre
over $T$, it is ACM and non-split.  Let $\mathcal{P}'=\coprod T$, the
union taken over all possible choices of $r>2$, $e$ and the
$a_i$'s. Thus, we get rank two non-split ACM bundles on all fibres of
$\mathcal{X}_{\mathcal{P}'}\to \mathcal{P}'$.

To prove that $\mathcal{P}'\to\mathcal{P}$ is dominant, it suffices to
show that the image of this map contains a dense set. Let
$x\in\mathcal{P}$ be a point such that $\mathcal{X}_x$ supports a rank
two non-split ACM bundle, say $E$. Let $r,e, a_i$ be the corresponding
invariants. Then by the universal property of Proposition
\ref{spreading1}, we see that there exists a point $y\in
S'=S'(r,e,a_i)$ such that $q(y)=x$ and the pull back of the
corresponding $G$ is this bundle $E$. By minimality of our resolution,
we see that $y\in S''$. Since $E$ is ACM, we see that $y\in T$. This
completes the proof.
\end{proof}

For a hypersurface $X\subset\mathbb{P}^n$ of degree $d$, the total
infinitesimal deformation $X_A$ of $X$ in $\mathbb{P}^n$ is contained
in $\mathbb{P}^n\times \mathrm{Spec}\, A$ where $A=k\oplus
V^{\vee}{}$, $V=\HH^0(\OX(d))$ and $V^{\vee}{}^2=0$ in the ring
$A$. The following corollary uses characteristic zero.

\begin{cor}\label{infinitesimal}
Assume that we have non-split rank two ACM bundles on a general
hypersurface of degree $d$ in $\mathbb{P}^n$ with $n\geq 4$. Then for
a general hypersurface $X$ of degree $d$ there exists a rank two
bundle $\mathcal{E}$ on $X_A $ such that $\mathcal{E}_{\mid X}$ is a
non-split ACM bundle.
\end{cor}
\begin{proof}
By theorem \ref{spreading}, under the hypothesis, we have a finite
type scheme $\mathcal{P}'$ mapping dominantly to $\mathcal{P}$ and a
rank two bundle on the `universal hypersurface' in $\mathbb{P}^n\times
\mathcal{P}'$ which is non-split and ACM on each fiber. Since
$\mathcal{P}$ is integral, by replacing $\mathcal{P}'$, we may clearly
assume that it is integral. Since the map is dominant, we may take a
generic multi-section and thus assume that
$\mathcal{P}'\to\mathcal{P}$ is generically finite. Replacing again
$\mathcal{P}'$ by a suitable open set, we may assume that the map is
etale since the characteristic is zero . If $x\in {\mathcal P}$
corresponding to $X\subset\mathbb{P}^n$ is the image of
$y\in\mathcal{P}'$, then the tangent spaces at $x$ and $y$ are
isomorphic. Now noting that $V$ is the tangent space at $x$ to
${\mathcal P}$, we are done.
\end{proof}

For the next theorem, we will need the following elementary lemma from
homological algebra, whose proof (which we omit) follows from the
construction of push-outs.
\begin{lemma}\label{pushout}
Let 
\[
\begin{array}{ccccc}
0&\to &\mathcal{A}&\stackrel{i}{\to}&\mathcal{B}\\
&&\downarrow j&&\downarrow j'\\
0&\to&\mathcal{A}'&\stackrel{i'}{\to}&\mathcal{B}'

\end{array}
\]
be a push-out diagram where
$\mathcal{A},\mathcal{A}',\mathcal{B},\mathcal{B}'$ are sheaves. Then
$i'$ splits if and only if there exists a homomorphism
$\alpha:\mathcal{B}\to \mathcal{A}'$ such that $\alpha\circ i=j$
\end{lemma}

Before we state the next theorem, let us fix some notation. Let $X$
be a smooth hypersurface of degree $d$ in $\mathbb{P}^n$ with $n\geq
3$ defined by $F=0$ and $E$ an ACM bundle on $X$ of arbitrary rank.


Let $G\in\HH^0(\OP(d))$ be such that the image $g\in V=\HH^0(\OX(d))$
is non-zero. Let $B=k[\epsilon]$ with $\epsilon^2=0$. The map
$-g:V^{\vee}{}\to k\epsilon$ defines a natural surjective ring
homomorphism $A\to B$ where $A$ is as above. Let $X_{\epsilon}\subset
\mathbb{P}^n\times \mathrm{Spec}\,B$ be the corresponding hypersurface
defined by $F-\epsilon G=0$. In the theorem we will assume that there
is a bundle $\mathcal{E}$ on $X_{\epsilon}$ such that the restriction
of $\mathcal{E}$ to $X\subset X_{\epsilon}$ is $E$.

As always, we have a resolution,
$$0\to F_1\sta{\phi} F_0\to E\to 0,$$ with the $F_i$'s direct sums of
line bundles on $\mathbb{P}^n$. This gives as before, by restriction
to $X$, a long exact sequence, and by splitting it to short exact
sequences, an exact sequence,
$$0\to E(-d)\to \ol{F_1}\to G\to 0$$ for some bundle $G$ on $X$ as in
(\ref{basic2}). Let us denote by $\zeta$ the corresponding extension
class in
$$\mathrm{Ext}^1_{\OX}(G,E(-d))=\HH^1(G^{\vee}\otimes E(-d)).$$

\begin{thm}\label{vanishing}
With the notation and assumptions as above, under the natural map
$$\HH^1(G^{\vee}\otimes E(-d))\stackrel{g}{\to} \HH^1(G^{\vee}\otimes
E))$$ the image of $\zeta$ is zero.
 
\end{thm}

\begin{proof}
Since $\mathcal{E}$ is flat over $k[\epsilon]$, we have,
$$E\cong \epsilon\mathcal{E}\cong \mathcal{E}/\epsilon\mathcal{E}.$$ 
 We get an exact sequence
$$0\to E=\epsilon\mathcal{E}\sta{i} \mathcal{E}\sta{\pi}
E=\mathcal{E}/\epsilon\mathcal{E}\to 0.$$ 

Let $p:X_{\epsilon}\to\mathbb{P}^n$ be the projection, which is
clearly a finite map. Taking the direct image under $p$  of the above
exact sequence and noting that $p$ restricted to $X\subset
X_{\epsilon}$ is a closed embedding in $\mathbb{P}^n$, we get

\begin{equation}\label{basic9} 
0\to E\stackrel{i}{\to} p_*\mathcal{E}=\mathcal{F}\stackrel{\pi}{\to} E\to 0.
\end{equation}
We want to interpret multiplication by $F$ on $\mathcal{F}$.  We know
that $F-\epsilon G$ annihilates $\mathcal{E}$, and so multiplication
by $F$ is the same as multiplication by $\epsilon G$. But
multiplication by $\epsilon$ is just the composite $i\circ\pi$. Thus
we see that multiplication by $F$ is the composite
\begin{equation}\label{basic10}
\mathcal{F}\stackrel{\pi}{\to} E\stackrel{g}{\to}
E(d)\stackrel{i}{\to}\mathcal{F}(d),
\end{equation} 
using the fact that $i\circ g=g\circ i:E\to \mathcal{F}(d)$.

The exact sequence (\ref{basic9}) gives an element $\eta$ in
$\mathrm{Ext}^1_{\OP}(E,E)$. From the
exact sequence
$$0\to F_1\to F_0\to E\to 0,$$ we get
$$\HH^0(F_1^{\vee}{}\otimes E)\to \mathrm{Ext}^1_{\OP}(E,E)\to
\HH^1(F_0^{\vee}{}\otimes E).$$ Since $E$ is ACM and $n\geq 3$ the
last term is zero. Thus we can lift $\eta$ to an $\alpha\in
\HH^0(F_1^{\vee}\otimes E)$ and we get the diagram
\[
\begin{array}{ccccccccc}
0&\to& F_1&\sta{\phi}& F_0&\to& E&\to& 0\\
&&\downarrow\alpha&&\downarrow\beta&&\|&&\\
0&\to &E&\sta{i}& \mathcal{F}&\sta{\pi}&E&\to&0
\end{array}
\]
Multiplication by $F$ gives the following commutative diagram 
\[
\begin{array}{ccc}
F_0(-d)&\sta{F}&F_0\\
\downarrow\beta&&\downarrow\beta\\
\mathcal{F}(-d)&\sta{F}&\mathcal{F}
\end{array}
\]
Since $F$ annihilates $E$, the top row factors as
$$F_0(-d)\sta{\psi} F_1\sta{\phi} F_0$$ for a suitable map
$\psi$. Using (\ref{basic10}) we get the following diagram, which we
will show is in fact commutative.

\[
\begin{array}{ccccccc}
F_0(-d)&&\sta{\psi}&&F_1&\sta{\phi}&F_0\\
\downarrow \beta&&&&\downarrow\alpha&&\downarrow\beta\\
\mathcal{F}(-d)&\sta{\pi}&E(-d)&\sta{g}&E&\sta{i}&\mathcal{F}
\end{array}
\]

We have
$$ig\pi\beta=F\beta=\beta F=\beta \phi\psi=i\alpha \psi.$$
Since $i$ is an inclusion, this implies, $g\pi\beta=\alpha
\psi$, proving commutativity. Restricting this diagram to $X$, we get

\[
\begin{array}{ccccccc}
\ol{F_0}(-d)&&\sta{\ol{\psi}}&&\ol{F_1}&\sta{\ol{\phi}}&\ol{F_0}\\
\downarrow\beta&&&&\downarrow\alpha&&\downarrow\beta\\
\mathcal{F}\otimes_{\OP}\OX(-d)&\sta{\pi}&E(-d)&\sta{g}&E&\sta{i}&
\mathcal{F}\otimes_{\OP}\OX 
\end{array}
\]

The image of $\ol{\psi}$ is just $E(-d)$ and thus we get a commutative
diagram

\[
\begin{array}{ccc}
E(-d)&\to&\ol{F_1}\\
\parallel &&\downarrow\alpha\\
E(-d)&\sta{g}&E
\end{array}
\]

The extension class $\zeta$ corresponds to the top row of the
following commutative diagram 
\[
\begin{array}{ccccccccc}
0&\to& E(-d)&\to&\ol{F_1}&\to&G&\to&0\\
&&\searrow g&&\swarrow\alpha&&&&\\
&&&E&&&&&
\end{array}
\]
and by lemma \ref{pushout}  under the natural map
$$\mathrm{Ext}^1(G,E(-d))\stackrel{g}{\to} \mathrm{Ext}^1(G,E),$$ $\zeta$
goes to zero.
Therefore under the map
$$\HH^1(G\,^{\vee}{}\otimes E(-d))\stackrel{g}{\to}
\HH^1(G\,^{\vee}{}\otimes E)$$
$\zeta$ goes to zero.  
\end{proof}

\begin{cor}\label{vanishing1}
Assume the rank of $E$ is two, $n\geq 4$ and that as before, $E$ can
be deformed in the direction of $g$. Then the natural map
$$\HH^2(E^{\vee}\otimes E(-d))\stackrel{g}{\to} \HH^2(E^{\vee}\otimes
E)$$ is zero.
\end{cor}

\begin{proof}
From  theorem \ref{vanishing} under the map,
$$\HH^1(G\,^{\vee}{}\otimes E(-d))\stackrel{g}{\to}
\HH^1(G\,^{\vee}{}\otimes E),$$ 
$\zeta$ goes to zero. By (\ref{basic2}), we have
$G\,^{\vee}{}=G(d-e)$. Thus under the map
$$\HH^1(G\otimes E(-e))\stackrel{g}{\to} \HH^1(G\otimes E(d-e)),$$
$\zeta$ goes to zero. If $E$ is indecomposable, $\HH^1(G\otimes
E(-e))$ is one dimensional with $\zeta$ as basis element by corollary
\ref{monogenic}. From (\ref{basic6}) we have,
$$\HH^1(G\otimes E(-e+k))=\HH^2(E\,^{\vee}{}\otimes E(-d+k)).$$
Thus we get that the map
$$\HH^2(E\,^{\vee}{}\otimes E(-d))\stackrel{g}{\to}
\HH^2(E\,^{\vee}{}\otimes E)$$ is zero.
\end{proof}
\end{section}

\section{Proof of theorem \ref{main}}

\begin{proof}[Proof of Theorem \ref{main} (1)]
For $n\geq 6$,
$$\HH^2(E\, ^{\vee}{}\otimes E(-d))=\HH^3(E\,^{\vee}{}\otimes G(-d))$$
by (\ref{basic7}). This group is zero by (\ref{basic8}). The proof now
follows from lemma \ref{suff1}.
\end{proof}

\begin{proof}[Proof of Theorem \ref{main} (2) a)]
For $n=5$, using (\ref{basic7}) with $i=1$, we have
$$\HH^1(X,E^{\vee}\otimes E(k))= \HH^2(X,E^{\vee}\otimes G(k)).$$ This
group is zero by (\ref{basic8}).
\end{proof}

\begin{proof}[Proof of Theorem \ref{main} (2) b)]
The proof is by contradiction. Assume that a general hypersurface $X$
of degree $d\geq 3$ has an indecomposable ACM bundle of rank two.
From corollaries \ref{infinitesimal} and \ref{vanishing1}, we see
that,
$$\HH^2(E\,^{\vee}{}\otimes E(-d))\stackrel{g}{\to}
\HH^2(E\,^{\vee}{}\otimes E)$$ is zero for any $g\in
\HH^0(\OX(d))$. By corollary \ref{monogenic} and (\ref{basic6}), we
know that the graded module $N=\oplus_k \HH^2(E\,^{\vee}{}\otimes
E(k))$ is cyclic and generated by a non-zero element $\zeta$ in degree
$-d$. Thus $N_0$ consists of multiples of $\zeta$ by elements $g\in
\HH^0(\OX(d))$. Since these are zero we get $N_i=0$ for $i\geq 0$. As
$n=5$, we have by Serre duality
$$\HH^2(E\,^{\vee}{}\otimes E(-d))\cong \HH^2(E\,^{\vee}{}\otimes
E(2d-6))^{\vee}.$$ 
Hence $N_{2d-6}\neq 0$.  If $d\geq 3$, then $2d-6\geq 0$ 
and this contradiction proves the result.
\end{proof}

\begin{proof}[Proof of Theorem \ref{main} (3)]
Assume that a general hypersurface $X$ of degree $d\geq 5$ has an
(indecomposable) ACM bundle of rank two. By the same arguments as in the
proof above, the graded module $N$ has $N_i=0$ for $i\geq 0$ and thus
$\HH^2(E\,^{\vee}{}\otimes E)=0$. By Serre duality,
we have 
$\HH^1(E^{\vee}\otimes E)\cong \HH^2(E\,^{\vee}{}\otimes E(d-5))=N_{d-5}$
and hence is zero for $d\geq 5$.

\end{proof}

\end{document}